# APPROXIMATION OF ANALYTIC
# FUNCTIONS WITH PRESCRIBED BOUNDARY
# CONDITIONS BY CIRCLE PACKING MAPS

Tomasz Dubejko

ABSTRACT. We use recent advances in circle packing theory to develop a constructive method for the approximation of an analytic function $F : \Omega \to \mathbf{C}$ by circle packing maps providing we have only been given $\Omega$, $|F'|\big|_\Omega$, and the set of critical points of $F$. This extends the earlier result of [CR] for $F$ with no critical points.

## 1 INTRODUCTION

In [CR] an inverse problem for circle packing and locally univalent function was considered. It was shown that for a given bounded simply connected domain $\Omega$ and a positive continuous function $\lambda : \partial\Omega \to (0, \infty)$ one can construct a sequence of locally univalent hexagonal circle packings and associated with it a sequence $\{f_n\}$ of circle packing maps such that $f_n \to f$ uniformly on compacta of $\Omega$ as $n \to \infty$, where $f$ is the unique (up to some standard normalization) locally univalent analytic function in $\Omega$ with $|f'| : \overline{\Omega} \to (0, \infty)$ being continuous and $|f'|\big|_{\partial\Omega} \equiv \lambda$.

Here we are interested in a generalization of the above problem, that is in approximation of an analytic function $F : \Omega \to \mathbf{C}$ via circle packing maps, where $F$ has a finite set of critical points in $\Omega$ and $|F'| : \overline{\Omega} \to [0, \infty)$ is continuous with $|F'| \equiv \lambda$ on $\partial\Omega$. We will prove that this can be achieved, roughly speaking, by taking the circle packing map from the portion of a regular hexagonal circle packing of suitably small mesh that fills up $\Omega$ to a combinatorially equivalent branched circle packing whose branch set approximates the set of critical points of $F$ and whose boundary circles have their radii determined by $\lambda$ and the mesh. The precise setup is laid out in Section 3 and the main result is contained in Theorem 3.1.

Techniques we employ here to verify our approximation scheme are ideologically quite different from the ones used in [CR]; the latter are closely linked with hexagonal combinatorics and with local univalence of functions considered, the former are based on advances in the theory of circle packing. In fact, arguments presented in this paper can be easily extended to other then hexagonal patterns (cf. [HR],[St]).

As the proof of Theorem 3.1 will heavily rely on results from the theory of circle packing, we combine them, together with an introduction of basic terminology, in

---

1991 *Mathematics Subject Classification.* 30G25, 30E25, 52C15.

The author gratefully acknowledges support of the Tennessee Science Alliance and the National Science Foundation. Research at MSRI is supported in part by grant no.DMS-9022140.

Section 2. The reader who would like to obtain more information about the subject should see [BSt] or [D1].

The author would like to acknowledge that both figures in this paper have been created with the help of `CirclePack`, a software package developed by Ken Stephenson.

## 2 Circle Packing Overview

Let $\mathbb{K}$ be a simplicial 2-complex which is a triangulation of a simply connected domain in the complex plane $\mathbf{C}$ and which has its orientation induced from the orientation of the plane. Write $\mathbb{K}^{(0)}$, int $\mathbb{K}^{(0)}$, bd $\mathbb{K}^{(0)}$, $\mathbb{K}^{(1)}$, and $\mathbb{K}^{(2)}$ for the sets of vertices, interior vertices, boundary vertices, edges, and faces of $\mathbb{K}$, respectively. A collection $\mathcal{P}$ of circles in $\mathbf{C}$ is said to be a *circle packing* for $\mathbb{K}$ if there is 1-to-1 correspondence between the elements of $\mathbb{K}^{(0)}$ and $\mathcal{P}$ ($\mathbb{K}^{(0)} \ni v \leftrightarrow C(v) \in \mathcal{P}$) such that if $\langle u,v \rangle \in \mathbb{K}^{(1)}$ then $C(u)$ and $C(v)$ are externally tangent and if $\langle u, v, w \rangle \in \mathbb{K}^{(2)}$ is a positively oriented triple then so is $\langle C(u), C(v), C(w) \rangle$.

Let $S_{\mathcal{P}}$ be a simplicial map $S_{\mathcal{P}} : \mathbb{K} \to \mathbf{C}$ which is first defined on $\mathbb{K}^{(0)}$ by mapping $v \in \mathbb{K}^{(0)}$ to the euclidean center of $C(v) \in \mathcal{P}$ and then is extended to $\mathbb{K}^{(1)}$ and $\mathbb{K}^{(2)}$ via barycentric coordinates. If $\triangle \in \mathbb{K}^{(2)}$ and $v$ is a vertex of $\triangle$ then $\alpha_{\mathcal{P}}(v, \triangle)$ will denote the angle at the vertex $S_{\mathcal{P}}(v)$ in the euclidean triangle $S_{\mathcal{P}}(\triangle)$. We define the *carrier* of $\mathcal{P}$ as carr$(\mathcal{P}) := S_{\mathcal{P}}(\mathbb{K})$ and the *angle sum* of $\mathcal{P}$ at a vertex $v$ as $\Theta_{\mathcal{P}}(v) := \sum_{\triangle \in \mathbb{K}^{(2)}} \alpha_{\mathcal{P}}(v, \triangle)$. It follows from our definition of a circle packing that $\Theta_{\mathcal{P}}(v)$ is a positive integer multiple of $2\pi$ if $v$ is an interior vertex of $\mathbb{K}$. Moreover, if $\mathcal{P}$ is *univalent*, i.e. the circles of $\mathcal{P}$ have mutually disjoint interiors, then $\Theta_{\mathcal{P}}(v) = 2\pi$ for each $v \in \text{int}\,\mathbb{K}^{(0)}$. However, the converse is not true (see Fig 1 (b)). If $v \in \text{int}\,\mathbb{K}^{(0)}$ and $\Theta_{\mathcal{P}}(v) = 2\pi n$, $n \geq 2$, then $v$ will be called a *branch vertex (point)* of $\mathcal{P}$ of *order $(n-1)$ (multiplicity $n$)*; geometrically this means that the circles in $\mathcal{P}$ associated with the neighboring vertices of $v$ wrap $n$-times around $C(v)$. The listing of all the branch vertices of $\mathcal{P}$ together with their orders will be called the *branch set* of $\mathcal{P}$ and denoted br$(\mathcal{P})$. A circle packing without branch points will be called *locally univalent*.

So far we have talked about circle packings but we have not said anything about their existence. Before we can state necessary and sufficient conditions for the existence we need the following definition, in which $\mathbf{N}$ will denote the set of non-negative integers.

**Definition 2.1.** *Let $\mathbb{K}$ be a finite or infinite triangulation of a simply connected domain in $\mathbf{C}$. A set $\{(v_1, l_1), (v_2, l_2), \dots\} \subset (\text{int}\,\mathbb{K}^{(0)}) \times \mathbf{N}$, possibly finite, is called a branch structure for $\mathbb{K}$ if every simple closed edge-path $\Gamma$ in $\mathbb{K}$ has at least $2\ell(\Gamma) + 3$ edges, where $\ell(\Gamma) = \sum_{i=1}^{n} \delta_i(\Gamma) l_i$ and $\delta_i(\Gamma)$ is equal to 1 if $v_i$ is enclosed by $\Gamma$ and 0 otherwise.*

The following theorem answers the existence question for circle packings (see [D1],[D2]).

**Theorem 2.2.** *Suppose $\mathbb{K}$ is a triangulation of a simply connected domain in $\mathbf{C}$. Let $b_1, \dots, b_m$ be interior vertices of $\mathbb{K}$, and let $k_1, \dots, k_m$ be non-negative integers. Then there exists a circle packing for $\mathbb{K}$ with branch set $\mathfrak{B} = \{(b_1, k_1), \dots, (b_m, k_m)\}$ if and only if $\mathfrak{B}$ is a branch structure for $\mathbb{K}$.*

(a)

(c)

(b)

Figure 1.   Different packings of the same complex:
(a) univalent, (b) locally univalent, (c) branched (and
its decomposition into univalent sheets)

Notice that when $0 = k_1 = \ldots = k_m$ then $\mathfrak{B}$ is a trivial case of a branch structure; in fact the number of pairs $(b_i, k_i)$ in $\mathfrak{B}$ is irrelevant in this case. Therefore we identify the sets $\big\{(b_1, 0), (b_2, 0), \ldots\big\}$ with the empty set.

For the purposes of this paper we will need a stronger version of Theorem 2.2 when $\mathbb{K}$ is finite which will be stated shortly, but first we have to introduce a function which is quite handy when working with circle packings. If $\mathcal{P}$ is a circle packing for $\mathbb{K}$ then a function $r_{\mathcal{P}} : \mathbb{K}^{(0)} \to (0, \infty)$ defined by $r_{\mathcal{P}}(v) :=$ the euclidean radius of $C(v)$ is called the *radius function* of $\mathcal{P}$. Now a stronger version of Theorem 2.2 (see [D1]).

**Theorem 2.3.** *Suppose $\mathbb{K}$ is a finite triangulation of a simply connected domain in $\mathbf{C}$. Let $b_1, \ldots, b_m$ be interior vertices of $\mathbb{K}$, and let $k_1, \ldots, k_m$ be non-negative integers. For any given $\rho : \mathrm{bd}\,\mathbb{K}^{(0)} \to (0, \infty)$ there exists a circle packing $\mathcal{P}$ for $\mathbb{K}$ with branch set $\mathfrak{B} = \big\{(b_1, k_1), \ldots, (b_m, k_m)\big\}$ and with $r_{\mathcal{P}}\big|_{\mathrm{bd}\,\mathbb{K}^{(0)}} \equiv \rho$ if and only if $\mathfrak{B}$ is a branch structure for $\mathbb{K}$. Moreover, $\mathcal{P}$ is unique up to isometries of $\mathbf{C}$.*

Suppose that $\mathcal{P}$ and $\mathcal{Q}$ are circle packings for $\mathbb{K}$. Then $\mathcal{P}$ and $\mathcal{Q}$ will be called *combinatorially equivalent* (or shortly c-equivalent) with complex $\mathbb{K}$. The map $F_{\mathcal{P},\mathcal{Q}} := S_{\mathcal{Q}} \circ S_{\mathcal{P}}^{-1} : \mathrm{carr}(\mathcal{P}) \to \mathrm{carr}(\mathcal{Q})$ will be called the *circle packing map* from $\mathcal{P}$ to $\mathcal{Q}$ (shortly, cp-map). The map $F_{\mathcal{P},\mathcal{Q}}^{\#} : \mathrm{carr}(\mathcal{P}) \to (0, \infty)$ defined on the set of vertices of $S_{\mathcal{P}}(\mathbb{K})$ by $F_{\mathcal{P},\mathcal{Q}}^{\#}(S_{\mathcal{P}}(v)) = \frac{r_{\mathcal{Q}}(v)}{r_{\mathcal{P}}(v)}$ and then extended affinely to faces of $S_{\mathcal{P}}(\mathbb{K})$ will be called the *ratio function* from $\mathcal{P}$ to $\mathcal{Q}$. The ratio function allows one to measure "distortion" of circle packings or how distinct two packings are. This is due to the following rather elementary fact, which gives, for example, the uniqueness in Theorem 2.3.

**Fact 2.4.**

(1) *If $\mathcal{P}$ and $\mathcal{Q}$ are circle packings for $\mathbb{K}$ with identical branch sets (possibly empty) then $F_{\mathcal{P},\mathcal{Q}}^{\#}$ does not obtain its infimum or supremum in $\mathrm{int}\,\mathbb{K}^{(0)}$ unless it is constant.*

(2) *If $\mathcal{P}$ and $\mathcal{Q}$ are circle packings for $\mathbb{K}$ such that $\mathrm{br}(\mathcal{P}) = \emptyset$ and $\mathrm{br}(\mathcal{Q}) \neq \emptyset$ then $F_{\mathcal{P},\mathcal{Q}}^{\#}$ does not attain its supremum in $\mathrm{int}\,\mathbb{K}^{(0)}$.*

We will now recall some results regarding circle packings and approximation of analytic functions. Let $\mathbb{H}$ be the regular hexagonal 2-complex of mesh 2 with vertices at $2k+l(1+\sqrt{3}i)$, $k, l \in \mathbf{Z}$. Write $\mathbb{H}_n := \frac{1}{n}\mathbb{H}$. Let $\mathcal{P}_n$ be the univalent circle packing whose carrier is $\mathbb{H}_n$, i.e. the regular hexagonal circle packing of circles of radius $\frac{1}{n}$. We notice that if $\mathfrak{B} = \{(b_1, k_1), \ldots, (b_m, k_m)\}$ is a branch structure for $\mathbb{H}_n$, $n \geq 1$, then necessarily $k_1 = \cdots = k_m = 1$, i.e. the $b_i$'s are simple branch points. Thus, when dealing with complexes $\mathbb{H}_n$, we adopt a convention by writing $\{b_1, \ldots, b_m\}$ for $\{(b_1, 1), \ldots, (b_m, 1)\}$. We now recall the definition of a discrete complex polynomial.

**Definition 2.5.** *A map $f : \mathbf{C} \to \mathbf{C}$ is a discrete complex polynomial for $\mathbb{H}_n$, $n \geq 1$, with the branch set $\mathfrak{B} = \{b_1, \ldots, b_m\}$ if the following are satisfied:*

(i) *there exists a circle packing $\mathcal{B}$ for $\mathbb{H}_n$ with the branch set $\mathfrak{B}$ such that $f$ is the cp-map from $\mathcal{P}_n$ to $\mathcal{B}$,*

(ii) *$f$ has a decomposition $f = \varphi \circ h$, where $h$ is a self-homeomorphism of $\mathbf{C}$ and $\varphi$ is a complex polynomial.*

*If $f$ is a discrete complex polynomial for $\mathbb{H}_n$ and $\mathcal{B}$ is as in (i) then $\mathcal{B}$ will be called the range packing of $f$.*

The following theorem is a result of Corollary 4.9 and Lemma 5.2 from [D2].

**Theorem 2.6.**

(1) *If $\mathfrak{B} = \{b_1, \ldots, b_m\}$ is a branch structure for $\mathbb{H}_n$, $n \geq 1$, then there exists a circle packing $\mathcal{B}$ for $\mathbb{H}_n$ with the branch set $\mathfrak{B}$ such that the cp-map $f : \mathbf{C} \to \mathbf{C}$ from $\mathcal{P}_n$ to $\mathcal{B}$ is a discrete complex polynomial.*

(2) *There exists a constant $\kappa \geq 1$ depending only on $m$ such that if $\mathcal{B}$ is a circle packing for $\mathbb{H}_n$, $n \geq 1$, and the cp-map from $\mathcal{P}_n$ to $\mathcal{B}$ has valence at most $m+1$ then $\frac{r_{\mathcal{B}}(w)}{r_{\mathcal{B}}(v)} \in \left(\frac{1}{\kappa}, \kappa\right)$ for any neighboring vertices $w, v$, where $r_{\mathcal{B}}$ is the radius function of $\mathcal{B}$. In particular, any discrete complex polynomial for $\mathbb{H}_n$, $n \geq 1$, with the branch set containing at most $m$ points is $K$-quasiregular, $K = K(m)$.*

The above theorem is a key factor in a proof of an approximation result for discrete complex polynomials which will be stated shortly. However, first we need to lay down some foundations. Let $F$ be a classical complex polynomial with the critical points $x_1, \ldots, x_m$ of orders $k_1, \ldots, k_m$, respectively, i.e. $\mathfrak{B} = \{(x_1, k_1), \ldots, (x_m, k_m)\}$ is the branch set of $F$. For each sufficiently large $n$ let $\mathfrak{B}_n = \{b_{1_1}(n), \ldots, b_{1_{k_1}}(n), \cdots, b_{m_1}(n), \ldots, b_{m_{k_m}}(n)\}$ be a branch structure for $\mathbb{H}_n$. Further, suppose the sequence $\{\mathfrak{B}_n\}$ has the property that $\lim_{n \to \infty} b_{i_j}(n) = x_i$ for $j = 1, \ldots, k_i$ and each $i$, $1 \leq i \leq m$. It is easy to see that such a sequence $\{\mathfrak{B}_n\}$ can be constructed for all sufficiently large $n$. As a consequence of Theorem 5.3 and Remark 5.4 of [D2] we have

**Theorem 2.7.** *For each sufficiently large $n$ there exists a discrete complex polynomial $f_n$ for $\mathbb{H}_n$ with the branch set $\mathfrak{B}_n$ such that the functions $f_n$ and $f_n^\#$ converge uniformly on compacta of $\mathbf{C}$ to $F$ and $|F'|$, respectively.*

The uniform convergence of $f_n^\# \to |F'|$ in the last theorem follows from the following result that we will need later (Theorem 1 of [DSt]).

**Theorem 2.8.** *Let $\Omega \subset \mathbf{C}$ be a bounded simple connected domain. Let $\{\mathcal{Q}_n\}$ and $\{\widetilde{\mathcal{Q}}_n\}$ be sequences of circle packings such that for each $n$ the packing $\mathcal{Q}_n$ is univalent, $\mathrm{carr}(\mathcal{Q}_n) \subseteq \Omega$, and $\mathcal{Q}_n$ and $\widetilde{\mathcal{Q}}_n$ are c-equivalent. Moreover, suppose that the sets $\mathrm{carr}(\mathcal{Q}_n)$ exhaust $\Omega$ and that the supremum of radii of circles in $\mathcal{Q}_n$ goes to 0 uniformly on compacta of $\Omega$ as $n \to \infty$. In addition, assume that functions $g_n$ converge uniformly on compacta of $\Omega$ to an analytic function $g : \Omega \to \mathbf{C}$ as $n \to \infty$, where $g_n : \mathrm{carr}(\mathcal{Q}_n) \to \mathrm{carr}(\widetilde{\mathcal{Q}}_n)$ is the cp-map from $\mathcal{Q}_n$ to $\widetilde{\mathcal{Q}}_n$. Then the ratio functions $g_n^\#$ converge uniformly on compacta of $\Omega$ to $|g'|$.*

## 3 The Main Result

In this section we will be concerned with the following problem. Let $\Omega \subset \mathbf{C}$ be a Jordan domain. Let $0 \in \Omega$ and $\xi \in \Omega$, $\xi > 0$. Suppose $\lambda : \partial\Omega \to (0, \infty)$ is a continuous function and $\mathfrak{B} = \{(x_1, k_1), \ldots, (x_m, k_m)\}$ is a subset of $\Omega \times \mathbf{N}$. It is a well-known fact that there exists the unique analytic function $F$ in $\Omega$ such that $F(0) = 0$, $F(\xi) > 0$, the set of critical points of $F$ is equal to $\mathfrak{B}$, and $|F'|$ has a continuous extension to $\overline{\Omega}$ with $|F'| \equiv \lambda$ on $\partial\Omega$. We are interested in the approximation of $F$ using cp-maps having only been given $\Omega$, $\lambda$, $\mathfrak{B}$, and $\xi$.

We will show that this can be achieved, roughly speaking, by taking the cp-map from a portion of a regular hexagonal circle packing of suitable small mesh that fills up $\Omega$ to a c-equivalent circle packing whose branch set approximates $\mathfrak{B}$ and whose radius function on the boundary is approximately equal to $\lambda$ times the mesh of the hexagonal packing. More precisely, write $\mathbb{O}_n$ for the maximal "complete" subcomplex of $\mathbb{H}_n$ contained in $\Omega$ (i.e., a simply connected simplicial 2-complex consisting of all faces of $\mathbb{H}_n$ whose closures are in $\Omega$, together with their edges and vertices). Denote $\mathcal{Q}_n$ the portion of $\mathcal{P}_n$ associated with $\mathbb{O}_n$, where $\mathcal{P}_n$ is as in Section 2. Suppose that for each sufficiently large $n$ there exists a circle packing $\widetilde{\mathcal{Q}}_n$ c-equivalent to $\mathcal{Q}_n$ such that

1. if $r_n : \mathbb{O}_n^{(0)} \to (0, \infty)$ is the radius function of $\widetilde{\mathcal{Q}}_n$, $v \in \mathrm{bd}\,\mathbb{O}_n^{(0)}$, and $z_v$ is a point on $\partial\Omega$ closest to $v$, then $r_n(v) = \frac{1}{n}\lambda(z_v)$,

2. if the set $\mathfrak{B}_n = \{b_{1_1}(n), \ldots, b_{1_{k_1}}(n), \cdots, b_{m_1}(n), \ldots, b_{m_{k_m}}(n)\}$ of distinct vertices of $\mathbb{O}_n$ is the branch set of $\widetilde{\mathcal{Q}}_n$ then the sequence $\{\mathfrak{B}_n\}$ has the property that $\lim_{n \to \infty} b_{i_j}(n) = x_i$ for $j = 1, \ldots, k_i$ and each $i$, $1 \leq i \leq m$,

(3) if $f_n : \text{carr}(\mathcal{Q}_n) \to \text{carr}(\widetilde{\mathcal{Q}}_n)$ is the cp-map from $\mathcal{Q}_n$ to $\widetilde{\mathcal{Q}}_n$ then $f_n(0) = 0$ and $f_n(\xi) > 0$,

then $\{\widetilde{\mathcal{Q}}_n\}$ will be called an *approximating sequence* of circle packings for $F$. We will prove the following theorem which is our main result.

**Theorem 3.1.** *Let $\Omega \subset \mathbf{C}$ be a Jordan domain. Let $0 \in \Omega$ and $\xi \in \Omega$, $\xi > 0$. Suppose $\lambda : \partial\Omega \to (0, \infty)$ is a continuous function and $\mathfrak{B} = \{(x_1, k_1), \ldots, (x_m, k_m)\}$ is a subset of $\Omega \times \mathbf{N}$. Let $F$ be the unique analytic function in $\Omega$ such that $F(0) = 0$, $F(\xi) > 0$, the set of critical points of $F$ is equal to $\mathfrak{B}$, and $|F'| : \overline{\Omega} \to [0, \infty)$ is continuous and $|F'|\big|_{\partial\Omega} \equiv \lambda$. Then there exists an approximating sequence $\{\widetilde{\mathcal{Q}}_n\}$ of circle packings for $F$. Moreover, if the maps $f_n$ are defined as in (3) above then $f_n$ and $f_n^{\#}$ converge uniformly on compacta of $\Omega$ to $F$ and $|F'|$, respectively.*

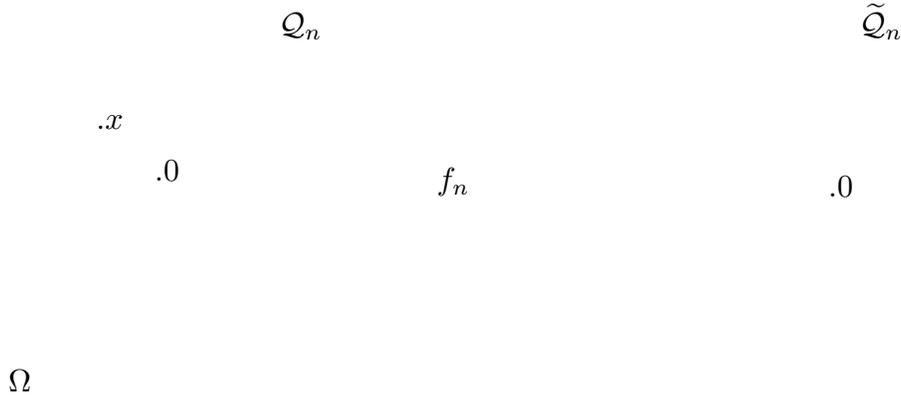

FIGURE 2. The circle packing map in an approximation sequence: $\lambda \equiv 1.4$, $\mathfrak{B} = \{(x, 1)\}$

Before we give a proof of the above theorem we want to make two remarks.

*Remark* 3.2. We observe that our construction of an approximating sequence $\{\widetilde{\mathcal{Q}}_n\}$ is based exclusively on $\Omega$, $\lambda$, $\mathfrak{B}$, and $\xi$, hence it gives a constructive method for the approximation of $F$.

*Remark* 3.3. If $\mathfrak{B} = \emptyset$ then we obtain the result of [CR].

To proof Theorem 3.1 we will need the following

**Lemma 3.4.** *Let $\Omega$ and $F$ be as in Theorem 3.1. Then there exists a sequence of polynomials $\{\chi_n\}$ such that $\chi_n \to F$ and $|\chi_n'| \to |F'|$ uniformly in $\overline{\Omega}$ as $n \to \infty$.*

*Proof.* First we need a formula for $F$. Let $u$ be the harmonic function in $\Omega$ with $u|_{\partial\Omega} = \log \lambda$. Denote by $v$ the harmonic conjugate of $u$ in $\Omega$. Write $\tau : \Omega \to D$ for the Riemann mapping with $\tau(0) = 0$ and $\tau(\xi) > 0$. Then the function $F$ is given by

$$F(z) = \int_0^z B(\tau(\eta)) e^{u(\eta) + iv(\eta)} d\eta,$$

where $B(z) = c \prod_{i=1}^m \left(\frac{z - \tau(x_i)}{1 - \overline{\tau(x_i)}z}\right)^{k_i}$ and $c$ is a suitable unimodular constant (such that $F(\xi) > 0$).

Let $\{\Omega_n\}$ be a sequence of Jordan domains such that $\overline{\Omega} \subset \Omega_{n+1} \subset \Omega_n$, $n \geq 1$, and the boundary of $\Omega_n$ converges to $\partial\Omega$ in the sense of Fréchet (cf. [LV, p.27],[H],[W]). Write $\tau_n : \Omega_n \to D$ for the Riemann mapping with $\tau_n(0) = 0$ and $\tau_n(\xi) > 0$. Then $\tau_n \to \tau$ uniformly in $\overline{\Omega}$. Define $F_n : \Omega_n \to \mathbf{C}$ by

$$F_n(z) = \int_0^z B(\tau_n(\eta))e^{u(\tau^{-1}\circ\tau_n(\eta))+iv(\tau^{-1}\circ\tau_n(\eta))}d\eta.$$

Then $F_n \to F$ uniformly in $\overline{\Omega}$. Moreover, as $|F_n'(z)| = |B(\tau_n(z))|e^{u(\tau^{-1}\circ\tau_n(\eta))}$ for $z \in \overline{\Omega}$, and $B$ and $u$ are continuous in $\overline{D}$ and $\overline{\Omega}$, respectively, we also get $|F_n'| \to |F'|$ uniformly in $\overline{\Omega}$. By Runge's theorem for each $n$ there exists a polynomial $\chi_n$ such that $\sup_{z\in\overline{\Omega}} |F_n(z) - \chi_n(z)| < \frac{1}{n}$ and $\sup_{z\in\overline{\Omega}} |F_n'(z) - \chi_n'(z)| < \frac{1}{n}$. The last implies the assertion of the lemma. $\square$

We are now ready to show Theorem 3.1.

*Proof of Theorem 3.1.* We first need to verify the existence of an approximating sequence of circle packings for the function $F$. This follows easily from Theorem 2.3, Definition 2.1, and the geometry of simplicial complexes $\mathbb{O}_n$. In fact, Definition 2.1 and the geometry of simplicial complexes $\mathbb{O}_n$ are only needed to obtain the condition (2) in the definition of an approximating sequence. The condition (3) of that definition is achieved immediately by applying translations and/or rotations, if required, to already constructed packings.

We will now show that if $\{\widetilde{Q}_n\}$ is an approximating sequence of circle packings for $F$ then the maps $f_n : \mathrm{carr}(Q_n) \to \mathrm{carr}(\widetilde{Q}_n)$ and $f_n^{\#} : \mathrm{carr}(Q_n) \to (0,\infty)$ converge uniformly on compacta of $\Omega$ to $F$ and $|F'|$, respectively. Let $G$ be a complex polynomial with the branch set $\mathfrak{B}$ such that $G(0) = 0$ and $G(\xi) = 1$. If $\mathfrak{B}_n$ is the branch set of $\widetilde{Q}_n$ then it is not hard to see that $\mathfrak{B}_n \subset \mathbb{O}_n^{(0)} \subset \mathbb{H}_n^{(0)}$ is a branch structure for $\mathbb{H}_n$ for all sufficiently large $n$. From Theorem 2.7 there exists, for each sufficiently large $n$, a discrete complex polynomial $g_n$ for $\mathbb{H}_n$ with the branch set $\mathfrak{B}_n$ such that $g_n \to G$ and $g_n^{\#} \to |G'|$ uniformly on compacta of $\mathbf{C}$ as $n \to \infty$. Write $R_n : \mathbb{H}_n^{(0)} \to (0,\infty)$ for the radius function of the range packing $\mathcal{U}_n$ of $g_n$. Since $|G'|\big|_{\partial\Omega} > 0$, $g_n^{\#}\big|_{\partial\mathbb{O}_n} \xrightarrow[n\to\infty]{} |G'|\big|_{\partial\Omega}$, and $g_n^{\#}(v) = R_n(v)/(\frac{1}{n})$, there exists $\sigma \geq 1$ such that $\frac{1}{n\sigma} \leq R_n\big|_{\mathrm{bd}\,\mathbb{O}_n^{(0)}} \leq \frac{\sigma}{n}$ for all $n$. Denote by $\widetilde{\mathcal{U}}_n$ the portion of $\mathcal{U}_n$ associated with the subcomplex $\mathbb{O}_n$ of $\mathbb{H}_n$. Recall that $r_n : \mathbb{O}_n^{(0)} \to (0,\infty)$ is the radius function of $\widetilde{Q}_n$. Since $\widetilde{\mathcal{U}}_n$ and $\widetilde{Q}_n$ are packings for $\mathbb{O}_n$ with the branch sets $\mathfrak{B}_n$, Fact 2.4 implies that $\frac{r_n}{R_n}(\cdot) : \mathbb{O}_n^{(0)} \to (0,\infty)$ attains its maximum and minimum in $\mathrm{bd}\,\mathbb{O}_n^{(0)}$. Hence

$$(3.1) \qquad \frac{1}{\sigma}\min_{\partial\Omega}\lambda \leq \frac{r_n}{R_n\big|_{\mathbb{O}_n^{(0)}}} \leq \sigma\max_{\partial\Omega}\lambda$$

because $\frac{1}{\sigma}\min_{\partial\Omega}\lambda \leq \frac{r_n}{R_n}\big|_{\mathrm{bd}\,\mathbb{O}_n^{(0)}} \leq \sigma\max_{\partial\Omega}\lambda$. Now it follows from Theorem 2.6 (applied to $g_n$) and (3.1) that there exists a constant $\kappa \geq 1$ such that

$$(3.2) \qquad \frac{1}{\kappa} \leq \frac{r_n(w)}{r_n(v)} \leq \kappa \quad \text{for any } n \text{ and any neighboring vertices } v, w \in \mathbb{O}_n^{(0)}.$$

The last conclusion and the construction of the maps $f_n$ imply that $\{f_n\}$ is a family of $K$-quasiregular mappings for some $K \geq 1$. Moreover, since $\Omega$ is bounded and by Fact 2.4 the ratio functions $f_n^\#$ are uniformly bounded by $\max_{\partial\Omega} \lambda$, it follows that the maps $f_n$ are uniformly bounded. Hence $\{f_n\}$ is a normal family. By taking a subsequence if necessary, assume that $f_n \to f$ uniformly on compacta of $\Omega$ as $n \to \infty$. We will first show that the limit function $f$ is not constant. If $f$ were constant then, by Theorem 2.8, functions $f_n^\#$ would converge uniformly on compacta of $\Omega$ to the constant zero-function. But this would contradict (3.1) for sufficiently large $n$ because $\left(\frac{r_n}{R_n}\right)\big|_{\mathbb{O}_n^{(0)}} = \left(\frac{f_n^\#}{g_n^\#}\right)\big|_{\mathbb{O}_n^{(0)}}$ and $|G'| > 0$ in $\Omega \setminus \{x_1, \ldots, x_m\}$. Hence $f$ is non-constant $K$-quasiregular mapping in $\Omega$.

We will now verify that $f$ is actually analytic. Since $f_n$ is bounded $K$-quasiregular mapping, by the Stoïlow's theorem ([LV]), it has a decomposition $f_n = \varphi_n \circ \psi_n$, where $\psi_n : \mathbb{O}_n \to D$ is a K-quasiconformal homeomorphism normalized by $\psi_n(0) = 0$ and $\psi_n(\xi) > 0$, and $\varphi_n : D \to \mathbf{C}$ is an analytic function with $\varphi_n(0) = 0$ and $\varphi_n(\psi_n(\xi)) > 0$. By extracting a subsequence from $\{f_n\}$ if necessary, we can assume that $\varphi_n \to \varphi$ and $\psi_n \to \psi$ uniformly on compacta of $\Omega$ and $D$, respectively, where $\varphi : \Omega \to D$ is a K-quasiconformal homeomorphism, $\psi : D \to \mathbf{C}$ is analytic, and $f = \varphi \circ \psi$. In particular, $\varphi$ and $\psi$ are not constant and the branch set of $\varphi$ is equal to $\big\{(\psi(x_1), k_1), \ldots, (\psi(x_m), k_m)\big\}$. Let $z_0 \in \Omega \setminus \{x_1, \ldots, x_m\}$. From the equicontinuity of the normal family $\{\psi_n\}$ and its uniform convergence on compacta to the function $\psi$ we get that there exist $\epsilon > 0$ and $n_0$ such that $f_n\big|_{B(z_0,\epsilon)}$ is 1-to-1 for $n > n_0$, where $B(z_0, \epsilon) = \{|z - z_0| < \epsilon\}$. Let $\widetilde{\mathcal{Q}}_n(z_0)$ be the portion of $\widetilde{\mathcal{Q}}_n$ associated with the largest complete hexagonal generation of $\mathbb{O}_n$ around $z_0$ contained in $B(z_0, \epsilon)$. Then $\{\widetilde{\mathcal{Q}}_n(z_0)\}$ is a sequence of univalent hexagonal circle packings with their number of generations around $z_0$ going to $\infty$ as $n \to \infty$. From the Hexagonal Packing Lemma of [RS] we conclude that the quasiconformal distortion of the maps $f_n$ at $z_0$ goes to 0 as $n \to \infty$. Hence the quasiconformal distortion of the mappings $\psi_n$ at $z_0$ goes to 0 as $n \to \infty$. Thus $\psi$ is 1-quasiconformal in $\Omega \setminus \{x_1, \ldots, x_m\}$. Since $\psi$ is a homeomorphism of $\Omega$, the last implies that $\psi$ is conformal in $\Omega$. Therefore $f = \varphi \circ \psi$ is analytic in $\Omega$ with the branch set $\mathfrak{B}$.

To complete our proof we need to show that $f = F$. To achieve this we have to prove that $|f'|$ has a continuous extension to $\partial\Omega$ and is equal to $\lambda$ there. We observe first that since $f$ is analytic in $\Omega$, by Theorem 2.8, $f_n^\# \to |f'|$ uniformly on compacta of $\Omega$ as $n \to \infty$. Let $\{\chi_n\}$ be a sequence of polynomials given by Lemma 3.4. Since $|F'|\big|_{\partial\Omega} > 0$ and $F$ has a finite branch set in $\Omega$, we can assume that for each $n$ the restriction of the branch set of $\chi_n$ to $\overline{\Omega}$ is equal to the branch set $\mathfrak{B}$ of $F$. Take $\delta > 0$ such that $|F'(z)| > 0$ for $z \in \Omega^\delta := \{z \in \overline{\Omega} : \mathrm{dist}(z, \partial\Omega) \leq \delta\}$. Given $\epsilon > 0$ let $N(\epsilon)$ be such that

$$(3.4) \qquad \begin{aligned} &\big||\chi_n'(z)| - |F'(z)|\big| < \epsilon, \quad z \in \overline{\Omega}, \quad n \geq N(\epsilon) \quad \text{and} \\ &\big|\tfrac{|\chi_n'(z)|}{|F'(z)|} - 1\big| < \epsilon, \quad z \in \Omega^\delta, \quad n \geq N(\epsilon). \end{aligned}$$

From the geometry of the complexes $\mathbb{H}_n$, Theorem 2.7, and the fact that $\mathfrak{B}_n$ is a branch structure for $\mathbb{O}_n$, it follows that there exists a sequence of discrete polynomials $\{p_n\}$ such that $p_n$ is a discrete polynomial for $\mathbb{H}_n$, $\mathrm{br}(p_n) \cap \overline{\Omega} = \mathfrak{B}_n$, and $p_n \to \chi_{N(\epsilon)}$ and $p_n^\# \to |\chi'_{N(\epsilon)}|$ uniformly on compacta of $\mathbf{C}$ as $n \to \infty$. Now (3.4) implies that there is $N'(\epsilon)$, $N'(\epsilon) \geq N(\epsilon)$, such that

$$(3.5) \qquad \left| \tfrac{p_n^\#}{|F'(z)|} - 1 \right| < 2\epsilon, \quad z \in \Omega^\delta, \quad n \geq N'(\epsilon).$$

Write $\mathcal{W}_n$ for the image packing of $p_n$. Denote by $\widetilde{\mathcal{W}}_n$ the portion of $\mathcal{W}_n$ associated with the subcomplex $\mathbb{O}_n$ of $\mathbb{H}_n$, and $\varrho_n$ the radius function of $\widetilde{\mathcal{W}}_n$. Since $\widetilde{\mathcal{W}}_n$ and $\widetilde{\mathcal{Q}}_n$ are circle packings for $\mathbb{O}_n$ with identical branch sets, Fact 2.4 shows that $\frac{\varrho_n}{r_n}(\cdot) : \mathbb{O}_n^{(0)} \to (0, \infty)$ has its maximum and minimum in bd $\mathbb{O}_n^{(0)}$. Hence $\left( \tfrac{p_n^\#}{f_n^\#} \right)\big|_{\mathbb{O}_n^{(0)}} : \mathbb{O}_n^{(0)} \to (0, \infty)$ has its maximum and minimum in bd $\mathbb{O}_n^{(0)}$. From the last conclusion, the construction of $\widetilde{\mathcal{Q}}_n$ (i.e., $f_n^\# \big|_{\partial(\mathrm{carr}(\mathcal{Q}_n))} \approx |F'| \big|_{\partial(\mathrm{carr}(\mathcal{Q}_n))} \approx \lambda$), the fact that $\mathrm{carr}(\mathcal{Q}_n) \subset \Omega^\delta$ for large $n$, and (3.5), it follows that there exists $N''(\epsilon)$, $N''(\epsilon) \geq N'(\epsilon)$, such that

$$(3.6) \qquad \left| \tfrac{p_n^\#}{f_n^\#}(v) - 1 \right| < 3\epsilon, \quad v \in \mathbb{O}_n^{(0)}, \quad n \geq N''(\epsilon).$$

We observe that $f_n^\# \big|_{\mathbb{O}_n^{(0)}} \leq \max_{\partial\Omega} \lambda$ by Fact 2.4(2). Hence, since $p_n^\#$ and $f_n^\#$ have been defined as simplicial extensions of $p_n^\# \big|_{\mathbb{O}_n^{(0)}}$ and $f_n^\# \big|_{\mathbb{O}_n^{(0)}}$, respectively, we obtain from (3.6)

$$(3.7) \qquad \left| p_n^\#(z) - f_n^\#(z) \right| < 3\epsilon \max_{\partial\Omega} \lambda, \quad z \in \mathrm{carr}(\mathcal{Q}_n), \quad n \geq N''(\epsilon).$$

Recall that $f_n^\# \to |f'|$ and $p_n^\# \to |\chi'_{N(\epsilon)}|$ uniformly on compacta of $\Omega$. Thus, by letting $n \to \infty$ in (3.7), one gets

$$\left| |\chi'_{N(\epsilon)}| - |f'(z)| \right| < 3\epsilon \max_{\partial\Omega} \lambda, \quad z \in \Omega.$$

As $\epsilon$ was arbitrary, (3.4) implies that $|F'| \equiv |f'|$ in $\Omega$. In particular, $|f'|$ has a continuous extension to $\partial\Omega$ and is equal to $|F'| \big|_{\partial\Omega} \equiv \lambda$ there. Hence $f \equiv F$ and the proof is complete. $\square$

*Concluding remarks.* 1) Notice that Theorem 3.1 gives, in essence, a method for the approximation of the integral of an analytic function on compacta of its domain.

2) Since the sequence of maps $\{f_n^\#\}$ converges uniformly on compacta of $\Omega$ to $|F'|$, $\{\log f_n^\#\}$ converges uniformly on compacta of $\Omega \setminus \{x_1, \ldots, x_m\}$ to the harmonic function $u$ in $\Omega \setminus \{x_1, \ldots, x_m\}$ which has isolated singularities at $x_i$ of degree $k_i$, $i = 1, \ldots, m$, and which satisfies boundary condition $u|_{\partial\Omega} = \log \lambda$. We also note that functions $f_n^\#$ are easier to construct then functions $f_n$; the latter require the radius functions of $\widetilde{\mathcal{Q}}_n$'s and the centers of circles in $\widetilde{\mathcal{Q}}_n$'s while the former require only radius functions of $\widetilde{\mathcal{Q}}_n$'s.

3) Arguments presented so far in this paper can easily be extended to other combinatorial patterns as follows. Let $\mathcal{O}$ be a univalent circle packing whose carrier is $\mathbf{C}$ (e.g. the "ball-bearing" packing – see [BSt] Figure 2(b)). Denote the geometric complex of $\mathcal{O}$ by $\mathbb{L}$. Assume that $\mathbb{L}$ is of bounded degree, i.e. there is a uniform bound on the number of neighboring vertices of each vertex in $\mathbb{L}$. Further, suppose $\mathcal{O}$ has the property that for each $\epsilon > 0$ there is $n_\epsilon$ such that all circles of $\frac{1}{n}\mathcal{O}$, $n \geq n_\epsilon$, contained in $\{|z| < 1\}$ have their radii at most $\epsilon$. Then, one can generalize Definition 2.5 and Theorem 2.6 to discrete complex polynomials for $\mathbb{L}_n$ $(= \frac{1}{n}\mathbb{L})$ as it was done in [D2]. Moreover, one can also obtain Theorem 2.7 for complexes $\mathbb{L}_n$ following the steps of the proof of Theorem 5.3 in [D2] and using Theorem 2.2 of [HR]. It is now a matter of replacing $\mathbb{H}_n$ by $\mathbb{L}_n$ to define approximation sequences based on complex $\mathbb{L}$ and to verify Theorem 3.1 for such sequences.

Department of Mathematics, University of Tennessee, Knoxville, TN 37996-1300
Mathematical Sciences Research Institute, Berkeley, CA 94720
*E-mail address*: tdubejko@msri.org